\documentclass{amsart}
\usepackage{amsmath,amsthm}
\usepackage{amsfonts,amssymb}
\usepackage{enumerate}

 \newtheorem{thm}{Theorem}[section]

\newtheorem{prop}[thm]{Proposition}

\newtheorem{defn}[thm]{Definition}

\theoremstyle{remark}

 \def\tr{{\triangle}}

\def\sph{\mathbb{S}^{d-1}}
\def\f{\frac}

 \def\a{{\alpha}}

 \def\k{{\kappa}}
 \def\t{{\theta}}
 \def\l{{\lambda}}
 
 \def\o{{\omega}}
 \def\s{{\sigma}}
 \def\la{{\langle}}
 \def\ra{{\rangle}}

 \def\CD{{\mathcal D}}
 
 \def\CH{{\mathcal H}}

 \def\CP{{\mathcal P}}
 
 \def\CT{{\mathcal T}}
 
 \def\CV{{\mathcal V}}
 
 \def\BB{{\mathbb B}}

 \def\NN{{\mathbb N}}

 \def\RR{{\mathbb R}}
  \def\SS{{\mathbb S}}
 \def\ZZ{{\mathbb Z}}
        
        \def\proj{\operatorname{proj}}
        
        \def\vi{\varphi}
        
  \def\ball{\mathbb{B}^{d}}

  \def\dist{\mathtt{d}}

\newcommand{\wt}{\widetilde}
\newcommand{\wh}{\widehat}

\makeindex             

\begin{document}

\title [Best Approximation on Spheres and Balls]
{Best Polynomial Approximation on the Unit Sphere and the Unit Ball}

\author{Yuan Xu}
\address{Department of Mathematics\\ University of Oregon\\
    Eugene, Oregon 97403-1222.}\email{yuan@uoregon.edu}

\maketitle

\begin{abstract}This is a survey on best polynomial approximation on the unit sphere and the unit ball. The central problem is 
to describe the approximation behavior of a function by polynomials via smoothness of the function. A major effort is to 
identify a correct gadget that characterizes smoothness of functions, either a modulus of smoothness or a $K$-
functional, the two of which are often equivalent. We will concentrate on characterization of best approximations, given 
in terms of direct and converse theorems, and report several moduli of smoothness and $K$-functionals, including recent 
results that give a fairly satisfactory characterization of best approximation by polynomials for functions in $L^p$ spaces, 
the space of continuous functions, and Sobolev spaces.
\end{abstract}
 
\section{Introduction}\label{sec:xu-intro}
\setcounter{equation}{0}

One of the central problems in approximation theory is to characterize the error of approximation of a function 
by the smoothness of the function. In this paper we give a short survey on best 
approximation by polynomials on the unit sphere $\sph$ and the unit ball $\BB^d$ in $\RR^d$ with 
$$
   \sph = \{x\in \RR^d: \|x\| =1\} \quad \hbox{and}\quad \BB^d = \{x: \|x\| \le 1\},
$$
where $\|x\|$ denotes the Euclidean norm of $x$. To get a sense of the main problem and its solution, let us consider
first  $\SS^1$ and $\BB^1$.

If we parametrize $\SS^1$ by $(\cos \t, \sin \t)$ with $\t \in [0, 2\pi)$ and identify a function $f$ defined on $\SS^1$ 
with the $2\pi$ periodic function $g(\t) = f(\cos \t, \sin \t)$, then polynomials on $\SS^1$ are precisely trigonometric
polynomials, so that polynomial approximation of functions on the circle $\SS^1$ is the same as trigonometric 
approximation of $2\pi$-periodic functions. Let $\CT_n$ denote the space of trigonometric polynomials of degree 
at most $n$,
$\CT_n: = \{a_0 + \sum_{k=1}^n a_k \cos k\t + b_k \sin k \t : a_k,b_k \in \RR\}$. Let $\|\cdot\|_p$ denote the 
$L^p(\SS^1)$ norm of $2\pi$-periodic functions on $[0,2\pi)$ 
if $1 \le p < \infty$, and the uniform norm of $C(\SS^1)$ 
if $p = \infty$. For $f \in L^p(\SS^1)$ if $1 \le p < \infty$, or $f \in C(\SS^1)$ if $p = \infty$, define 
$$
   E_n(f)_p : = \inf_{t_n \in \CT_n} \| f - t_n\|_p,
$$ 
the error of best approximation by  trigonometric polynomials. The convergence behavior of $E_n(f)_p$ is usually 
characterized by a modulus of smoothness. 
For  $f \in L^p(\SS^1)$ if $1 \le p < \infty$ or $f \in C(\SS^1)$ if $p = \infty$, $r = 1, 2, \ldots$ and $t > 0$, the
modulus of smoothness defined by the forward difference is 
$$
   \o_r (f;t)_p := \sup_{|\t| \le t}  \left \| \overrightarrow{\tr}_\t^r f \right \|_p, \qquad 1 \le p \le \infty,
$$
where $\overrightarrow{\tr}_h f(x):=  f(x+h) - f(x)$ and $\overrightarrow{\tr}_h^r := 
\overrightarrow{\tr}_h^{r-1} \overrightarrow{\tr}_h$. The characterization of best approximation on $\SS^1$ is 
classical (cf. \cite{DL, Timan}). 
 
\begin{thm}
For  $f \in L^p(\SS^1)$ if $1 \le p < \infty$ or $f \in C(\SS^1)$ if $p = \infty$,
 \begin{equation} \label{JacksonS1}
       E_n(f)_p \le c \, \o_r \left (f; n^{-1}\right )_p, \quad 1 \le p \le \infty, \quad n =1,2,....
\end{equation}
On the other hand,
\begin{equation} \label{converseS1}
       \o_r(f; n^{-1})_p \le c\, n^{-r} \sum_{k=1}^n k^{r-1} E_{k-1}(f)_p, \qquad 1 \le p \le \infty.
\end{equation}
\end{thm}

The theorem contains two parts. The direct inequality \eqref{JacksonS1} is called the Jackson estimate, 
its proof requires constructing a trigonometric polynomial that is close to the best approximation. The 
weak converse inequality \eqref{converseS1} is called the Bernstein estimate as its proof relies on 
the Bernstein inequality. Throughout this paper, we let $c$, $c_1$, $c_2$ denote constants independent
of $f$ and $n$.  Their values may differ at different times. 

Another important gadget, often easier to use in theoretical studies, is the $K$-functional defined by
$$
 K_r(f,t)_p := \inf_{g \in W_p^r} \left \{\|f - g\|_p + t^r \|g^{(r)} \|_{p} \right \},
$$
where $W_p^r$ denotes the Sobolev space of functions whose derivatives up to $r$-th order are all in $L^p(\SS^1)$. 
The modulus of smoothness $ \o_r (f,t)_p$ and the K-function $ K_r(f,t)_p$ are known to be equivalent:
for some constants $c_2 > c_1 > 0$, independent of $f$ and $t$, 
\begin{equation}\label{m-Kequiv}
    c_1 K_r (f,t)_p \le   \o_r (f,t)_p  \le c_2 K_r (f,t)_p.  
\end{equation}

All characterizations of best approximation, either on the sphere $\sph$ or on the ball $\ball$, encountered in
this paper follow along the same line: we need to define an appropriate modulus of smoothness and use it 
to establish direct and weak converse inequalities; and we can often define a $K$-functional that is 
equivalent to the modulus of smoothness. 

\medskip

{\bf Convention:} In most cases, our direct and weak converse estimates are of the same form as those in \eqref{JacksonS1} 
and \eqref{converseS1}. In those cases, we shall simply state that the direct and weak converse theorems hold and 
will not state them explicitly. 

\medskip

We now turn our attention to approximation by polynomials on the interval $\BB^1 := [-1,1]$. Let $\Pi_n$ 
denote the space of polynomials of degree $n$ and let $\|\cdot\|_p$ also denote the $L^p$ norm of 
functions on $[-1,1]$ as in the case of $\SS^1$. For $f \in L^p(\BB^1)$, $1 \le p < \infty$, or $f \in C(\BB^1)$
for $p = \infty$, define 
$$
   E_n(f)_p : = \inf_{t_n \in \Pi_n} \| f - p_n\|_p, \qquad 1 \le p \le \infty.  
$$
The difficulty in characterizing $E_n(f)_p$ lies in the difference between approximation behavior at the
interior and at the boundary of $\BB^1$. It is well known that polynomial approximation on $\BB^1$ 
displays a better convergence behavior at points close to the boundary than at points in the interior. 
A modulus of smoothness that is strong enough for both direct and converse estimates should catch this
boundary behavior. 

There are several successful definitions of modulus of smoothness in the literature. The most satisfactory 
one is due to Ditzian and Totik in \cite{Di-To}. For $r \in \NN$ and $h > 0$, let $\wh \tr_h^r$ denote the central 
difference of increment $h$, defined by 
\begin{equation} \label{central-diff}
   \wh \tr_h f(x)  = f(x+\tfrac{h}2) - f(x- \tfrac{h}2) \quad \hbox{and} \quad 
           \wh \tr_h^r  = \wh \tr_h^{r-1} \tr, \quad r =2, 3, \ldots.
\end{equation}
Let $\varphi(x) := \sqrt{1-x^2}$. For $r =1,2\ldots$, and $1\leq p\leq \infty$, the Ditzian-Totik moduli of smoothness 
are defined  by
\begin{equation}\label{modu-DT}
   \o_\vi^r (f,t)_{p} := \sup_{0 < h \le t}  \left \| \wh \tr_{h \varphi}^r f \right \|_{L^p[-1,1]},
\end{equation}
where $\wh \Delta_{h\vi(x)}^r f(x) = 0$ if $x\pm {r h\vi(x)}/2\notin [-1,1]$. Both direct theorem and 
weak converse theorem 
for $E_n(f)_p$ hold for this modulus of smoothness. Furthermore, the $K$-functional that is equivalent 
to this modulus of smoothness is defined by, for $t > 0$ and $r =1,2\ldots$, 
\begin{equation}\label{K-func-DT}
     K_{r,\vi} (f,t)_{p} : = \inf_{g\in C^r[-1,1]}  \left \{ \|f-g\|_{p} +
         t^r \|\varphi^{r} g^{(r)}\|_{p} \right \}. 
\end{equation}
  
\medskip

In the rest of this paper, we discuss characterization of the best approximation on the sphere $\sph$
and on the ball $\BB^d$. The problem for higher dimension is much harder. For example,  
functions on $\sph$ are no longer periodic, and there are interactions between variables for functions 
on $\sph$ and $\BB^d$. 

The paper is organized as follows. The characterization of best approximation on the sphere is 
discussed in the next section, and the characterization on the ball is given in Section 3. In Section 4 
we discuss recent result on Sobolev approximation on the ball, which are useful for spectral methods
for numerical solution of partial differential equations. The paper ends with a problem on characterizing 
best polynomial approximation of functions in Sobolev spaces. 

\section{Approximation on the Unit Sphere}
\label{sec:xu_approx_sphere}

We start with necessary definitions on polynomial spaces and differential operators. 

\subsection{Spherical harmonics and spherical polynomials} 
\label{subsec:xu-SphericaHarmonics}

For $\sph$ with $d \ge 3$, spherical harmonics play the role of trigonometric functions for the unit circle. 
There are many books on spherical harmonic -- we follow \cite{DaiXu13}. 
Let $\CP_n^d$ denote the space of real  homogeneous polynomials of degree $n$ and let $\Pi_n^d$ 
denote the space of real polynomials of degree at most $n$. It is known that 
\begin{equation*}
  \dim \CP_n^d = \binom{n+d-1}{n} \quad\hbox{and} \quad \dim \Pi_n^d = \binom{n+d}{n}.
\end{equation*}
Let $\Delta := \partial_1^2 + \cdots + \partial_d^2$ denote the usual Laplace operator. A polynomial
$P \in \Pi_n^d$ is called harmonic if $\Delta P = 0$. For $n=0,1,2,\ldots $ let  $\CH_{n}^{d}:
=\left\{ P\in \CP_{n}^{d}:\Delta P=0\right\}$ be the linear space of real harmonic polynomials that are
homogeneous of degree $n$. Spherical harmonics are the restrictions of elements in $\CH_n^d$ on
the unit sphere. It is known that 
$$
a_n^d: = \dim \CH_n^d = \dim \CP_n^d - \dim \CP_{n-2}^d.
$$
Let $\Pi_n^d(\sph)$ denote the space of polynomials restricted on $\sph$. Then 
$$
   \Pi_n^d(\sph) = \bigoplus_{0 \le j \le n/2} \CH_{n-2j}^d \Big \vert_{\sph} \quad \hbox{and} \quad
   \dim \Pi_n^d(\sph) = \dim \CP_{n}^{d}+\dim \CP_{n-1}^{d}. 
$$
For $x \in \RR^d$, write $x = r \xi$, $r \ge 0$, $\xi \in \sph$. The Laplace operator can be written as 
\begin{equation*} 
\Delta = \frac{\partial^2}{\partial r^2} + \frac{d - 1}{r} \frac{\partial}{\partial r} + \frac{1}{r^2}\Delta_0,
\end{equation*}
where $\Delta_0$ is a differential operator on $\xi$, called the Laplace-Beltrami operator; see
\cite[Section 1.4]{DaiXu13}. The spherical harmonics are eigenfunctions of $\Delta_0$. More precisely,
$$
   \Delta_0 Y(\xi) = -n (n+d-2) Y (\xi), \qquad  Y \in \CH_n^d.
$$

The spherical harmonics are orthogonal polynomials on the sphere. Let  $d \s$ be the
surface measure, and 
$\o_{d-1}$ be the surface area of $\sph$. For $f, g \in L^1(\sph)$, define
$$
    \la f, g \ra_{\sph}  := \f{1}{\o_{d-1}} \int_{\sph} f(\xi) g(\xi) d\s(\xi).
$$ 
If 
$Y_n \in \CH_n^d$ for $n = 0,1,\ldots$, then $\la Y_n, Y_m \ra_{\sph} =0$ if $n \ne m$. 
A basis $\{Y_\nu^n: 1 \le \nu \le a_n^d\}$ of $\CH_n^d$ is called orthonormal if 
$\la Y_\nu, Y_\mu \ra_{\sph} = \delta_{\nu,\mu}$. In terms of an orthonormal basis, the
reproducing kernel $Z_{n,d}(\cdot, \cdot)$ of $\CH_n^d$ can be written as 
$Z_{n,d}(x,y) = \sum_{1 \le \nu \le a_n^d} Y_\nu(x) Y_\nu(y)$, and the addition formula for the
spherical harmonics states that 
\begin{equation}\label{zonal}
   Z_{n,d}(x,y) =  \frac{n+\l }{\l} C_n^\l (\la x,y\ra), \qquad \l = \f {d-2} 2,
\end{equation}
where $C_n^\l$ is the Gegenbauer polynomial of one variable. If $f \in L^2(\sph)$, then the Fourier
orthogonal expansion of $f$ can be written as 
$$
  f = \sum_{n=0}^\infty \proj_n f, \qquad \proj_n: L^2(\sph) \mapsto \CH_n^d,
$$
where the projection operator $\proj_n$ can be written as an integral 
$$
  \proj_n f(x) = \f{1}{\o_{d-1} } \int_{\sph} f(y) Z_{n,d}(x,y) d\s(y). 
$$

For $f \in L^p(\sph)$, $1 \le p < \infty$, or $f \in C(\sph)$ if $p = \infty$, the error of best approximation
by polynomials of degree at most $n$ on $\sph$ is defined by 
$$
   E_n(f)_p := \inf_{P \in \Pi_n(\sph)} \| f- P\|_{p}, \qquad 1 \le p \le \infty, 
$$
where the norm $\|\cdot\|_p$ denote the usual $L^p$ norm on the sphere and 
$\|\cdot\|_\infty$ denote the uniform norm on the sphere. Our goal is to characterize this quantity in terms 
of some modulus of smoothness. The direct theorem of 
such a characterization requires a polynomial that is close to the least polynomial that approximates $f$. 
For $p =2$, the $n$-th polynomial of best approximation is the partial sum, 
$$
    S_n f = \sum_{k=0}^n \proj_k f,
$$ 
of the Fourier orthogonal expansion, as the standard Hilbert space theory shows. For $p \ne 2$, a polynomial
of near best approximation can be given in terms of a cut--off function, which is a $C^\infty$-function $\eta$ 
on $[0,\infty)$ such that $\eta(t)=1$ for $0\le \eta(t)\le 1$ 
and $\eta(t)=0$ for $t \ge 2$. If $\eta$ is such a function, define
\begin{equation}\label{Vnf}
 S_{n,\eta} f(x):= \sum_{k=0}^{\infty} \eta\left(\f k{n}\right) \proj_k f(x). 
\end{equation}
Since $\eta$ is supported on $[0,2]$, the summation in $ S_{n,\eta} f$ can be terminated at $k =2 n-1$,
so that $ S_{n,\eta} f$ is a polynomial of degree at most $2n-1$. 

\begin{thm} \label{thm:near-best}
Let $f \in L^p(\sph)$ if $1 \le p < \infty$ and $f \in C(\sph)$ if $p = \infty$. Then
\begin{enumerate}[  \quad \rm(1)]
 \item $ S_{n,\eta}  f \in \Pi_{n}(\sph)$ and $ S_{n,\eta} f = f $ for $f \in \Pi_n^d(\sph)$.
 \item For $n \in \NN$, $\| S_{n,\eta}  f\|_{p} \le c \|f\|_{p}$.
 \item For $n \in \NN$, there is a constant $c > 0$, independent of $f$, such that 
\begin{equation*} 
        \|f -  S_{n,\eta}  f\|_{p} \le (1+c) E_{n}(f)_{p}.
\end{equation*}
\end{enumerate}
\end{thm}
This near-best approximation was used for approximation on the sphere already in \cite{Kam} and it 
has become a standard tool by now. For further information, including a sharp estimate of its kernel
function, see \cite{DaiXu13}. 
 
\subsection{First Modulus of Smoothness and $K$-functional}
\label{subsec:xu-2first}

The first modulus of smoothness is defined in terms of spherical means. 

\begin{defn}\label{1st-modulusS}
For $0 \le \t \le \pi$ and $f \in L^1(\sph)$, define the spherical means 
\begin{equation*}
   T_\theta f(x) := \frac{1}{\o_{d-1}} \int_{\SS_x^\bot}f(x \cos \theta + u\sin\theta ) d\s(u),
\end{equation*}
where $\SS_x^\bot := \{y \in \SS^{d-1}: \langle x,y\rangle =0\}$.
 For $f\in L^p(\sph)$, $1 \le p < \infty$, or $C(\sph)$, $p = \infty$,
and $ r > 0$, define 
\begin{equation}\label{omega*}
  \o_r^*(f, t)_p := \sup_{|\t| \le t} \| (I - T_\t)^{r/2} f\|_p,
\end{equation}
where $(I-T_\t)^{ r/2}$ is defined by its formal infinite series when $ r/2$ is not an integer.
\end{defn}
 
The equivalent $K$-functional of this modulus is defined by 
\begin{equation}\label{K-func*}
  K_r^*(f,t)_p := \inf_g \left\{ \|f-g\|_p + t^r  \left \|(- \Delta_0)^{r/2}g \right\|_p \right\},
\end{equation}
where $\Delta_0$ is the Laplace-Beltrami operator on the sphere and the infimum is taken over all
$g$ for which $(-\Delta_0)^{r/2} g \in L^p(\sph)$.

This modulus of smoothness was first defined and studied in \cite{BBP, P}.

\begin{thm} \label{thm:approxS1}
For $1 \le p \le \infty$, the modulus of smoothness $\o_r^*(f,t)_p$ can be used to establish both direct and 
weak converse theorems, and it is equivalent to $K_r^*(f,t)_p$.
\end{thm}

The direct and the weak converse theorems were established in various stages by several authors 
(see \cite{BBP, Kal, NL1, P, WL} and \cite{NL, WL} for further references), before it was finally established 
in full generality by Rustamov \cite{Rus1}. A complete proof is given in \cite{WL} and a simplified proof
can be found in \cite{DaiXu13}. 

The spherical means $T_\t$ are multiplier operators of Fourier orthogonal series,
i.e.,
\begin{equation}\label{eq:multiplier}
   \proj_n T_\t f = \frac{C_n^{\l}(\cos\t)}{C_n^{\l}(1)} \proj_n f, \qquad \l = \frac{d-2}{2}, \quad n=0,1,2,\ldots. 
\end{equation}
This fact plays an essential role for studying this modulus of smoothness.  

It should be mentioned that this multiplier approach can be extended to weighted approximation 
on the sphere, in which $d\s$ is replaced by $h_\k^2 d\s$ , where $h_\k$ is a function invariant 
under a reflection group. The simplest such weight function is of the form
$$
    h_\k(x) = \prod_{i=1}^d |x_i |^{\k_i}, \qquad \k_i \ge 0, \quad x \in \sph, 
$$
when the group is $\ZZ_2^d$. Such weight functions were first considered by Dunkl associated 
with Dunkl operators. An extensive theory of harmonic analysis for orthogonal expansions 
with respect to $h_\k^2 (x)d\s$ has been developed (cf. \cite{DaiXu10, DX}), in parallel with the classical 
theory for spherical harmonic expansions. The weighted best approximation in $L^p(h_\k^2; \sph)$
norm was studied in \cite{X05}, where analogues of the modulus of smoothness $\o_r^*(f,t)_p$
and K-functional $K_r^*(f,t)_p$ are defined with $\|\cdot\|_p$ replaced by the norm of $L^p(h_\k^2; \sph)$
for $h_\k$ invariant under a reflection group, and a complete analogue of Theorem \ref{thm:approxS1}
was established. 

\medskip

The advantages of the moduli of smoothness $\o_r^*(f,t)_p$ are that they are well--defined for all $r > 0$ 
and they have a relatively simple structure through multipliers. These moduli, however, are difficult 
to compute even for simple functions.  

\subsection{Second  Modulus of Smoothness and $K$-functional}
\label{subsec:xu-2second}

The second modulus of smoothness on the sphere is defined through rotations on the sphere. 
Let $SO(d)$ denote the group of orthogonal matrix of determinant 1. For $Q \in SO(d)$, 
let $T(Q) f(x) := f(Q^{-1}x)$. For $t > 0$, define 
$$
  O_t := \left \{Q \in SO(d):  \max_{x \in \SS^{d-1}} \dist(x,Qx) \le t \right \},
$$
where $\dist(x,y): = \arccos \la x, y\ra$ is the geodesic distance on $\SS^{d-1}$.

\begin{defn}
For $f\in L^p(\sph)$, $1 \le p < \infty$, or $C(\sph)$, $p = \infty$, and $r > 0$, define 
 \begin{equation}\label{omegaD}
\wt \o_r(f, t)_p : = \sup_{Q \in O_t} \| \tr_Q^r f\|_p,
    \qquad\hbox{where}\quad  \tr_Q^r: = (I - T_Q)^r.
\end{equation}
\end{defn}

For $r =1$ and $p=1$, this modulus of smoothness was introduced and used in \cite{CWZ} and further 
studied in \cite{KW}. For studying best approximation on the sphere, these moduli were introduced 
and investigated by Ditzian in \cite{Di1} and he defined them for more general spaces, including 
$L^p(\SS^{d-1})$ for $p>0$. 

\begin{thm} \label{thm:approxS2}
The modulus of smoothness $\wt \o_r(f,t)_p$ can be used to establish both direct and weak converse 
theorems for $1 \le p \le \infty$, and it is equivalent to the $K$-functional 
$K_r^*(f,t)_p$ for $1 < p < \infty$, 
but the equivalence fails if $p =1$ or $p=\infty$. 
\end{thm}

The direct and weak converse theorems were established in \cite{Di2} and \cite{Di1}, respectively. The 
equivalence of $\wt \o_r(f;t)_p$ and $ K_r^*(f,t)_p$ for $1< p < \infty$ was proved in \cite{DDH}, and the 
failure of the equivalence for $p = 1$ and $\infty$ was shown in \cite{Di3}. 

The equivalence passes to the moduli of smoothness and shows, in particular, that $\wt \o_r(f;t)_p$ is 
equivalent to the first modulus of smoothness $\o_r^*(f;t)_p $ for $1 < p < \infty$ but not for 
$p = 1$ and $p = \infty$.  

 \medskip
 
One advantage of the second moduli of smoothness $\wt \o_r(f;t)_p $ is that they are independent of the 
choice of coordinates. These moduli, however, are also difficult to compute even for fairly simple functions. 
 
\subsection{Third modulus of smoothness and $K$-functional}
\label{subsec:xu-2third}

The third modulus of smoothness on the sphere is defined in terms of moduli of 
smoothness of one variable on multiple circles. For $1\le i, j \le d$, we let $\tr_{i,j,t}^r$ 
be the $r$-rh forward difference acting on the angle of the polar coordinates on the 
$(x_i,x_j)$ plane. For instance,  take $(i,j) =(1,2)$ as an example,
\begin{equation*} 
\tr_{1,2, \t}^r f(x) =  \overrightarrow{\tr}_\t^r
   f \left (x_1 \cos(\cdot)- x_2 \sin (\cdot), x_1 \sin(\cdot) + x_2 \cos (\cdot),
        x_3,\ldots,x_d \right).
\end{equation*}
Notice that if $(x_i,x_j)= s_{i,j} (\cos \t_{i,j}, \sin \t_{i,j})$ then 
$$
(x_1 \cos \t - x_2 \sin \t, x_1 \sin \t + x_2 \cos \t) = s_{i,j} \cos (\t_{i,j} + \t),
$$
so that $\tr_{1,2, \t}^r f(x)$ can be regarded as a difference on the circle of the $(x_i,x_j)$ plane. 

\begin{defn} \label{def:modulus}
For $r  =1,2, \ldots$, $t > 0$, and $f \in L^p(\SS^{d-1})$, $1 \le p < \infty$, or $f \in C(\SS^{d-1})$ 
for $p = \infty$, define
\begin{equation} \label{eq:modulus}
 \o_r (f,t)_p : =  \max_{1 \le i < j \le d} \sup_{ |\t| \le t} \left \|\tr_{{i,j,\t}}^r f \right \|_p.
\end{equation}
\end{defn}

The equivalent $K$-functional is defined using the angular derivative  
\begin{equation*} 
        D_{i,j} : = x_i \partial_j - x_j \partial_i = \frac{\partial}{\partial \t_{i,j}}, \quad 1 \le i\ne j \le d
\end{equation*}
where $\t_{i,j}$ is the angle of polar coordinates in $(x_i,x_j)$-plane defined as above. For 
$r \in \NN_0$ and $t > 0$, the $K$-functional is defined by 
\begin{equation} \label{eq:K-func-sphere}
   K_r(f,t)_p : = \inf_{g } \left\{ \|f - g\|_p + t^r \max_{1 \le i<j \le d}
        \|D_{i,j}^r g\|_p\right\},
\end{equation}
where $g$ is taken over all $g  \in L^p(\sph)$  for which $D_{i,j}^r g \in L^p(\sph)$ for all $1\le  i,j \le d$. 

\begin{thm} \label{thm:approxS3}
The modulus of smoothness $\o_r(f,t)_p$ can be used to establish both direct and weak converse 
theorems, and is equivalent to $K_r(f,t)_p$ for $1 \le p \le \infty$.
\end{thm}

These moduli and $K$-functionals were introduced in \cite{DaiXu10}, where the above theorem was 
proved. Furthermore, it was also shown that 
$$
    K_r(f, n^{-1})_p \sim \|f - S_{n,\eta} f\|_p +
             n^{-r} \max_{1\le i<j \le d} \|D_{i,j}^r S_{n,\eta} f\|_p, 
$$
where $S_{n,\eta}$ is the polynomial defined in \eqref{Vnf}. 

For comparison with the other two moduli of smoothness, it was proved in \cite{DaiXu10} 
that
for $r =1,2, \ldots$ and $1 \le p \le \infty$, 
$$
  \o_r (f, t)_p \le \wt \o_r(f, t)_p, \qquad   0< t< 1. 
$$ 
Furthermore, for $1 < p < \infty$, the two moduli of smoothness are equivalent if $r =1$ or $r= 2$.
Thus, the direct theorem with $\o_r(f,t)_p$ is at least not weaker than the one with either one of the
other two moduli of smoothness. Furthermore, all three moduli are equivalent if $1 < p < \infty$ and 
$r =1$ or $2$. It remains an open problem if $\o_r(f,t)_p$ is  equivalent to other two moduli of 
smoothness for $1 < p < \infty$ and $r \ge 3$ or for $p =1$ and $p = \infty$. 

The angular derivatives are related to the Laplace-Beltrami operator by  
   \begin{equation*}
           \Delta_0 = \sum_{1 \le i<j \le d} D_{i,j}^2. 
\end{equation*}
Since the $K$-functional $K^*_r(f,t)_p$ is defined in terms of $\Delta_0$ and the $K$-function $K_r(f,t)$ is 
defined in terms of $D_{i,j}$, it indicates that $K_r(f,t)_p$ may be stronger than $K_r^*(f,t)_p$ if we
believe that the parts encode more information than the whole. 

The main advantage of the modulus of smoothness $ \o_r (f, t)_p$ lies in the fact that it is defined in
terms of moduli of smoothness of one variable, which allows us to tap into the well established theory
of trigonometric approximation of one variable, and it also means that $\o_r(f,t)_p$ can be computed
relatively easily (see \cite{DaiXu10} for examples). 

One interesting phenomenon observed from the computational example is that the best approximation 
on $\sph$ for $d\ge 3$ displays a boundary behavior rather like approximation by polynomials on $[-1,1]$. 
This is not all that surprising on second thought, but it does put $d=2$ in approximation on $\SS^{d-1}$ apart
from $d \ge 3$.  

\section{Approximation on the Unit Ball}
\label{sec:xu_approx_ball}
 
On the unit ball, we often work with weighted approximation with a fairly general weight function. We shall
restrict our discussion to the classical weight function 
$$
w_\mu(x): = (1-\|x\|^2)^{\mu -1/2}, \qquad \mu > -1/2, \quad x \in \ball, 
$$
for which the most has been done. We start with an account of orthogonal structure. 

\subsection{Orthogonal structure on the unit ball}
\label{subsec:xu-ortho-ball}

For the weight funciton $W_\mu$, we consider the space $L^p(w_\mu, \ball)$ for $1 \le p < \infty$ 
or $C(\BB^d)$ when $p = \infty$. The norm of the space  $L^p(w_\mu, \ball)$ will be denoted by $\|f\|_{\mu,p}$, 
taken with the measure $w_\mu(x) dx$. The inner product of $L^2(w_\mu, \ball)$ is defined by
$$
  \la f, g\ra_{\mu,p}: = b_\mu \int_{\ball} f(x) g(x) w_\mu(x) dx,
$$ 
where $b_\mu$ is the normalization constant of $w_\mu$ such that $\la 1, 1\ra_{\mu,p} =1$. 
Let $\CV_n^d(w_\mu)$ denote the space of polynomials of degree $n$ that are orthogonal to
polynomials in $\Pi_{n-1}^d$ with respect to the inner product $\la \cdot,\cdot \ra_{\mu,p}$. It is 
known that $\dim \CV_n^d(w_\mu)= \binom{n+d-1}{n}$.  The 
orthogonal polynomials in $\CV_n^d(w_\mu)$ are eigenfunctions 
of a second order differential operator:  for $g \in \CV_n^d(w_\mu)$, 
\begin{equation}\label{CDmu}
     \CD_{\mu} g: = \big(\Delta - \la x,  \nabla\ra^2 - (2 \mu+ d-1) \la x, \nabla \ra \big)g = 
       - n(n+2 \mu +d-1) g. 
\end{equation}

For $\nu \in \NN_0^d$ with $|\nu|= n$, let $P_\nu^n$ denote an orthogonal polynomial in $\CV_n^d(w_\mu)$. 
If $\{P_\nu^n: |\nu| =n\}$ is an orthonormal basis of $\CV_n^d$, then the reproducing kernel 
$P_n(w_\mu; \cdot, \cdot)$ of $\CV_n^d(w_\mu)$ can be written as
$P_n(w_\mu; x, y) = \sum_{|\nu|=n} P_\nu^n (x) P_\nu^n(y)$. This kernel satisfies a closed-form 
formula (\cite{X99}) that will be given later in this subsection. Let $L^2(w_\mu, \ball)$, then the Fourier 
orthogonal expansion of $f$ can be written as 
$$
  f = \sum_{n=0}^\infty \proj_n^\mu f, \qquad \proj_n^\mu: L^2(w_\mu, \ball) \mapsto \CV_n^d(w_\mu),
$$
where the projection operator $\proj_n$ can be written as an integral 
$$
  \proj_n^\mu f(x) = b_\mu \int_{\ball} f(y) P_n(w_\mu; x,y) w_\mu(y) dy. 
$$

For $f \in L^p(w_\mu, \BB^d)$, $1 \le p < \infty$, or $f \in C(\ball)$ if $p = \infty$, the error of best approximation
by polynomials of degree at most $n$ is defined by 
$$
   E_n(f)_{\mu,p} := \inf_{P \in \Pi_n^d} \| f- P\|_{\mu,p}, \qquad 1 \le p \le \infty. 
$$
The direct theorem for $E_n(f)_{\mu,p}$ is also established with the help of a polynomial that is a near best approximation to $f$. For $p =2$, the best polynomial of degree $n$ is again the partial sum,
$S_n^\mu f:= \sum_{k=0}^n \proj_k^\mu f$, of the Fourier orthogonal expansion, whereas for $p \ne 2$ we 
can choose the polynomial as 
\begin{equation}\label{Vnfball}
 S_{n,\eta}^\mu f(x):= \sum_{k=0}^{\infty} \eta\left(\f k{n}\right) \proj_k^\mu f(x), 
\end{equation}
where $\eta$ is a cut--off function as in \eqref{Vnf}. The analogue of Theorem \ref{thm:near-best} holds for
$S_{n,\eta}^\mu$ and $\|\cdot\|_{\mu,p}$ norm. 

If $\mu$ is an integer or a half integer, then the orthogonal structure of $L^2(w_\mu, \ball)$ is closely related to the 
orthogonal structure on the unit sphere, which allows us to deduce many properties for analysis on the unit ball
from the corresponding results on the unit sphere. The connection is based on the following identity: 
if $d$ and $m$ are positive integers, then for any $f\in L(\SS^{d+m-1})$, 
$$
\int_{\SS^{d+m-1}} f(y) d\s_{d+m} = \int_{\BB^d} (1-\|x\|^2)^{\frac{m-2}{2}}
 \biggl[\int_{\SS^{m-1}} f\Big(x, \sqrt{1-\|x\|^2} \xi\Big) d\s_m(\xi)\biggr] dx.
$$
This relation allows us to relate the space $\CV_n^d(w_\mu)$ with $\mu = \f{m-1}{2}$ directly to a subspace
of $\CH_n^{d+m}$, which leads to a relation between the reproducing kernels. 

For $\mu = \frac{m-1}{2}$, the reproducing kernel $P_n(w_\mu; \cdot,\cdot)$ satisfies, for $m > 1$,
\begin{align*} 
  P_n(w_\mu; x, y)  = \frac{1}{\o_m} \int_{\SS^{m-1}} Z_{n,d+m} \Big( (x,x'), (y, \sqrt{1-\|y\|^2} \xi)\Big)
         d\s_m(\xi),
\end{align*}
where $(x,x')\in \SS^{d+m-1}$ with $x \in \BB^d$ and $x' = \|x'\| \xi \in \BB^m$ with $\xi \in \SS^{m-1}$, 
and it satisfies, for $m=1$ and $y_{d+1} = \sqrt{1-\|y\|^2}$, 
$$
 P_n(w_0; x, y)   = \f12 \left[Z_{n,d+m}  \big ((x,x'), (y, y_{d+1})\big) +
     Z_{n,d+m} \big( (x,x'), (y, - y_{d+1}) \big)\right]. 
$$
Using the identity \eqref{zonal}, we can then obtain a closed-form formula for $P_n(w_\mu; \cdot, \cdot)$,
which turns out to hold for all real $\mu > -1/2$. 

\subsection{First Modulus of Smoothness and $K$-functional}
\label{subsec:xu-3first}

The first modulus of smoothness on the unit ball is an analogue of $\o_r^*(f,t)_p$ on the sphere, 
defined in the translation operator $T_\t^\mu$. Let $I$ denote the identity matrix and
$$
  A(x) := (1-\|x\|^2) I + x^T x, \quad x =(x_1,\ldots, x_d) \in \ball.
$$
For $W_\mu$ on $\ball$, the generalized translation operator is given by 
\begin{align*} 
  T_\theta^\mu f(x) =  b_\mu (1-\|x\|^2)^{\frac{d-1}{2}} &
     \int_{\Omega} f\big(\cos \theta x + \sin \theta \sqrt{1-\|x\|^2}\,u\big ) \left(1 - u A(x) u^T \right)^{\mu-1} du,
\end{align*}
where $\Omega$ is the ellipsoid $\Omega = \{u: u A(x) u^T \le 1\}$ in $\RR^d$. 

\begin{defn}\label{defn:modulus1-Bd}
Let $f\in L^p(W_\mu, \BB^d)$ if $1\leq p<\infty$, and $f\in C(\BB^d)$ if $p=\infty$. For $r =1,2,\ldots,$ and 
$t > 0$, define 
\begin{equation*}
  \o_r^* (f, t)_{\mu,p}:=\sup_{|\t|\leq t} \|\tr_{\t,\mu}^r f\|_{p,\k}, \qquad  \tr_{\t,\mu}^r f:=\big(I-T_\t^\mu \big)^{r/2} f.
\end{equation*}
\end{defn}

The equivalent $K$-functional is defined via the differential operator $\CD_\mu$ in \eqref{CDmu},
$$
  K_r^*( f, t)_{\mu,p} : = \inf_{g} \big \{ \|f-g\|_{\mu,p} + t^r \|\CD_\mu^r g\|_{\mu,p} \big\},
$$
where $g$ is taken over all $g\in L^p(W_\mu, \BB^d)$ for which  $\CD_\mu^r g\in L_p(W_\mu, \ball)$.

\begin{thm} \label{thm:approxB1}
For $1 \le p \le \infty$, the modulus of smoothness $\o_r^*(f,t)_{\mu,p}$ can be used to establish both
direct and weak converse theorems, and it is equivalent to $K_r^*(f,t)_{\mu,p}$.
\end{thm}

These moduli of smoothness and $K$-functionals were defined in \cite{X05} and Theorem \ref{thm:approxB1} 
was also proved there. The integral formula of $T_\t^\mu f$ was found in \cite{X05c}. In fact, these results 
were established for more general weight functions of 
$h_\k^2 w_\mu$ with $h_\k$ being a reflection invariant function. The operator $T_\theta^\mu$ is a 
multiplier operator and satisfies 
$$
\proj_n^\mu \big(T_\t^\mu f \big)  = \frac{C_n^{\l_\mu}(\cos \t)}{C_n^{\l_\mu}(1)}  \proj_n^\mu f, 
\qquad \l_\mu = \mu + \frac{d-1}2, \quad n =0,1,\ldots,
$$
which is an analogue of \eqref{eq:multiplier}. The proof of Theorem \ref{thm:approxB1} can be 
carried out following the proof of Theorem \ref{thm:approxS1}. 

The advantage of the moduli of smoothness $\o_r^\ast (f, t)$ are that they are well--defined for all $r > 0$ 
and their connection to multipliers, just like the first moduli of smoothness on the sphere. These moduli, 
however, are difficult to compute even for simple functions.  

\subsection{Second  Modulus of Smoothness and $K$-functional}
\label{subsec:xu-3second}

The second modulus of smoothness is inherited from the third moduli of smoothness on the sphere. 
With a slight abuse of notation, we write $w_\mu(x):=(1-\|x\|^2)^{\mu-\f12}$ for either the weight
function on $\BB^d$ or that on  $\BB^{d+1}$, and write $\tr^r_{i,j,\t}$ for either the difference
operator on $\RR^d$ or that on $\RR^{d+1}$. This should not cause any confusion from the
context. We denote by $\wt f$ the extension of $f$ defined by 
\begin{equation*} 
     \wt f(x,x_{d+1}) = f(x), \qquad (x,x_{d+1}) \in \BB^{d+1}, \quad x \in \BB^d.
\end{equation*}

\begin{defn} \label{defn:modulus-B}
Let $\mu = \f{m-1}{2}$, $f \in L^p(w_\mu, \BB^d)$ if $1 \le p < \infty$ and $f \in C(\BB^d)$ if $p = \infty$. 
For $r =1,2\ldots,$ and $t>0$, define 
\begin{align*} 
  \o_r(f,t)_{p,\mu} :=  \sup_{|\t|\le t} 
  \Big \{ \max_{1 \le i<j \le d}  \| \tr_{i,j,\t}^r  f\|_{L^p(\BB^{d},W_ {\mu})},
\max_{1 \le i \le d}  \| \tr_{i,d+1,\t}^r \wt f\|_{L^p(\BB^{d+1},W_ {\mu-1/2})} \Big\},
\end{align*}
where for $m=1$, $ \| \tr_{i,d+1,\t}^r \wt f\|_{L^p(\BB^{d+1},W_ {\mu-1/2})}$
is replaced by  $\| \tr_{i,d+1,\t}^r \wt f\|_{L^p(\SS^d)}$.
\end{defn}

The equivalent $K$-functional is defined in terms of the angular derivatives
$D_{i,j}$, and 
is defined for all $\mu \ge 0$ by 
\begin{align*}  
     K_r(f,t)_{p,\mu} := \inf_{g\in C^r(\BB^d)} \Big\{  & \|f-g\|_{L^p( W_\mu; \BB^d)}
    + t^r   \max_{1 \le i<j \le d}  \| D_{i,j}^r  g\|_{L^p(W_ {\mu}; \BB^d)}  \\
  & \qquad\qquad +  t^r \max_{1 \le i \le d}  \| D_{i,d+1}^r  \wt g\|_{L^p(W_ {\mu-1/2}; \BB^{d+1})} \Big \},
\notag
\end{align*}
where if $\mu = 0$, then $\| D_{i,d+1}^r  \wt g\|_{L^p(W_ {\mu-1/2}; \BB^{d+1})}$ is replaced by 
$\| D_{i,d+1}^r  \wt g\|_{L^p(\SS^d)}$.

\begin{thm} \label{thm:approxB2}
Let $\mu = \frac{m-1}{2}$. For $1 \le p \le \infty$, the modulus of smoothness $\o_r(f,t)_{\mu,p}$ can be 
used to establish both direct and weak converse theorems, and it is equivalent to $K_r(f,t)_{\mu,p}$. 
\end{thm}

The moduli of smoothness $\o_r(f,t)_{p,\mu}$ and the $K$-functionals $K_r(f,t)_{p,\mu}$ were 
introduced in \cite{DaiXu10} and Theorem \ref{thm:approxB2} was proved there. The proof relies 
heavily on the correspondence between $L^p(w_\mu, \BB^d)$ and $L^p(\SS^{d+m-1})$. In the 
definition of $\o_r(f,t)_{\mu,p}$, the term that involves the difference of $\wt f$ may look strange but
it is necessary, since $\tr_{i,j,\t}^r$ are differences in the spherical coordinates. 

For comparison with the first modulus of smoothness $\o_r^*(f,t)_{\mu,p}$, we only have that 
for $1<p<\infty$, $r = 1,2, \ldots$ and  $0<t<1$, 
\begin{equation*}
     \o_r(f,t)_{p,\mu} \leq c \o_r^\ast (f, t)_{p,\mu}.
\end{equation*}
In all other cases, equivalences are open problems.  Furthermore, the main results are established
only for $\mu = \f{m-1}{2}$, but they should hold for all $\mu \ge  0$ and perhaps even $\mu > -1/2$,
which, however, requires a different proof from that of \cite{DaiXu10}.

One interesting corollary is that, for $d =1$, $\o_r(f,t)_{\mu,p}$ defines a
modulus of smoothness on $\BB^1 = [-1,1]$ that is previously unknown. For $\mu = \frac{m-1}{2}$, this 
modulus is given by, for $f \in L^p(w_\mu, [-1,1])$,  
\begin{equation*}
\o_r(f,t)_{p,\mu} := \sup_{|\t| \le t}
  \left(c_\mu \int_{\BB^2} \left |{\tr}_\t ^r
    f(x_1 \cos (\cdot)  + x_2 \sin (\cdot)) \right|^p  w_{\mu-\f12}(x) dx \right)^{1/p}. 
\end{equation*}

One advantage of the moduli of smoothness is that they can be relatively easily computed. Indeed, they
can be computed just like the second modulus of smoothness on the sphere; see \cite{DaiXu10} for 
several examples.  

\subsection{Third modulus of smoothness and $K$-functional}  

The third modulus of smoothness on the unit ball is similar to $\o_r(f,t)_{p,\mu}$,
but with the term that
involves the difference of $\wt f$ replaced by another term that resembles the difference in the Ditzian--Totik modulus 
of smoothness. To avoid the complication of the weight function, we state this modulus of smoothness only
for $\mu =1/2$ for which $w_\mu(x) =1$. In this subsection, we write $\|\cdot\|_p : = \|\cdot \|_{1/2,p}$. 

Let $e_i$ be the $i$-th coordinate vector
of $\RR^d$ and let  $\wh\Delta_{h e_i}^r$ be the $r$-th central difference in the
direction of $e_i$. More precisely,
$$
\wh \Delta_{h e_i} f (x):= f(x+he_i) - f(x - h e_i), \qquad \wh \Delta^{r+1}_{he_i} f(x)  = \wh  \Delta_{he_i} \wh \Delta^{r}_{he_i} f(x).
$$
As in the case of $[-1,1]$, we assume that $\wh\Delta_{h e_i}^r$ is zero if either  of the points 
$x \pm r \frac{h}{2} e_i$ does not belong to $\BB^d$.  

\begin{defn} 
Let $f \in L^p(\BB^d)$ if $1 \le p < \infty$ and $f \in C(\BB^d)$
if $p = \infty$. For $r =1,2,\ldots$ and $t>0$,
\begin{equation*} 
   \o_\vi^r(f,t)_{p} := \sup_{0 < |h| \le t} \Big \{ \max_{1 \le i<j \le d}
      \| \tr_{i,j,h}^r  f\|_{p},
      \max_{1 \le i \le d}  \|\wh \tr_{h\varphi  e_i}^r f\|_{p} \Big \}.
\end{equation*}
\end{defn}

With $\varphi(x):=\sqrt{1-\|x\|^2}$, 
the equivalent $K$-functional is defined by
\begin{align*}  
   K_{r,\vi}(f,t)_{p}
   := \inf_{g\in W_p^r(\BB^d)} \Big\{ \|f-g\|_{p}
        + t^r   \max_{1 \le i<j \le d}  \| D_{i,j}^r  g\|_{p} +
      t^r \max_{1 \le i \le d}  \| \varphi^r \partial_i^r g\|_{p} \Big \}.
\end{align*}

\begin{thm} \label{thm:approxB3}
For $1 \le p \le \infty$, the modulus of smoothness $\o_\vi^r (f,t)_{\mu,p}$ can be used to establish both 
direct and weak converse theorems, where the direct estimate takes the form
\begin{equation*} 
  E_n (f)_{p} \le c\, \o_\vi^r(f,n^{-1})_p+n^{-r}\|f\|_{p} 
\end{equation*}
in which the additional term $n^{-r }\|f\|_p$ can be dropped when $r= 1$, and it is equivalent to $K_{r,\vi}(f,t)$ in 
the sense that 
$$
      c^{-1}\o_\vi^r(f,t)_p\le   K_{r,\vi}(f,t)_p\leq c\, \o_\vi^r(f,t)_p+c\,t^r\|f\|_p,
$$
where the term $t^r\|f\|_p$ on the right side can be dropped when $r=1$.
\end{thm}

These moduli of smoothness and $K$-functionals were also defined in \cite{DaiXu10}, and e
Theorem \ref{thm:approxB3} was proved there. For $d =1$, they agree with the Ditzian--Totik moduli 
of smoothness and $K$-functionals. The $K$-functional $ K_{r,\vi}(f,t)_{\mu,p}$ can be defined by
replacing $\|\cdot\|_p$ with $\|\cdot\|_{\mu,p}$ in the definition of $K_{r,\vi}(f,t)_{p}$, which were 
used to prove direct and weak converse theorems for $E_n(f)_{\mu,p}$ in terms of the $K$-functionals
in \cite{DaiXu10}. 

For comparison with the second $K$-functional $K_r(f,t)_{\mu,p}$, which is only defined for 
$\mu = \frac{m-1}{2}$, $m = 1,2,\ldots$, we know that for $1 \le p \le \infty$, 
\begin{align*} 
    K_{1,\vi}(f,t)_{\mu,p} \sim K_1(f,t)_{\mu,p}
\end{align*}
and, for $r > 1$, there is a $t_r > 0$ such that
\begin{align*} 
       K_r(f,t)_{\mu,p} \le c \, K_{r,\vi}(f,t)_{\mu,p} + c \, t^r \|f\|_{\mu,p},  \quad 0<t<t_r,
\end{align*}
where we need to assume that $r $ is odd if $p = infty$. We can also state the result for comparison 
of the moduli of smoothness $\o_{r,\vi}(f,t)_{p}$ and $\o_r(f,t)_{1/2,p}$ accordingly. The other direction of
the equivalence for $r =2, 3, \ldots$ remain open. 

The advantages of the modulus of smoothness $\o_\vi^r(f,t)_p$ and the $K$-functional $\o_\vi^r(f,t)_p$
are that they are more intuitive,  as direct extensions of the Ditizian--Totik modulus of smoothness
and $K$-functional, and that the modulus of smoothness is relatively easy to compute.

\section{Approximation in the Sobolev Space on the Unit Ball}
\label{sec:xu_approx_Sobolev}

For $r= 1,2,\ldots$ we consider the Sobolev space $W_r^p(\ball)$ with the norm defined by
$$
   \|f\|_{W_r^p(\ball)} =  \Big (\sum_{|\a| \le r} \|\partial^\a f\|_p \Big)^{1/p}. 
$$
The direct theorem given in terms of the $K$-functional yields immediately an estimate of $E_n(f)_p$ for
functions in the Sobolev space. In the spectral method for solving partial differential equations, we often
want estimates for the errors of derivative approximation as well. In this section, we again 
let $\|\cdot\|_p = \|\cdot \|_{1/2,p}$.

Approximation in Sobolev space requires estimates of derivatives. One such result was proved in
\cite{DaiXu11}, which includes the following estimates 
\begin{align*}
 \|D_{i,j}^r (f - S_n^\mu f)\|_{p,\mu} & \le c  E_n ( D_{i,j}^r f)_{p,\mu}, \quad 1 \le i < j \le d,
\end{align*}
and a similar estimate that involves $D_{i,d+1} \wt f$. However, what we need is an estimate that 
involves only derivatives $\partial^\a$ instead of $D_{i,j}^r$. In this regard, the following result
can be established.

\begin{prop}
If $f \in W_p^{s}(\ball)$ for $1 \le p < \infty$, or $f \in C^{s}(\ball)$ for $p = \infty$,  then for $|\a| = {s}$, 
\begin{align} \label{eq:approx-weight}
  \big \| \phi^{|\a|/p} (\partial^\a f - \partial^\a S_{n,\eta} f) \big \|_{p} \le c E_{n-|\a|} (\partial^\a f)_{p}
    \le c n^{-s} \|f\|_{W_p^s(\ball)}, 
\end{align}
where $S_{n,\eta} f = S_{n,\eta}^{1/2} f$ is the near--best approximation defined in \eqref{Vnfball}.  
\end{prop}

The estimate \eqref{eq:approx-weight} in the proposition, however, is still weaker than what is needed in 
the spectral method, which requires an estimate similar to \eqref{eq:approx-weight} but without the term
$[\phi(x)]^{|\a|/p}= (1-\|x\|^2)^{|\a|/p}$. It turns out that the near--best approximation $S_{n,\eta} $ is 
inadequate for obtaining such an estimate. What we need  is the
orthogonal structure of the Sobolev space $W_2^r(\ball)$. 

The orthogonal structure of  $W_2^r(\ball)$ was studied first in \cite{X08} for the
case $r=1$, and in \cite{PX,X06} for the case $r=2$, and in \cite{LX} for general $r$. The inner 
product of $W_2^r(\ball)$ is defined by 
$$
\la f, g\ra_{-s}:= \la \nabla^s f, \nabla^s g \ra_{\ball} + \sum_{k=0}^{\lceil \f s 2 \rceil -1} \la \Delta^k f, \Delta^k g \ra_{\sph}. 
$$
Let $\CV_n^d(w_{-s})$ denote the space of polynomials of degree $n$ that are orthogonal to
polynomials in $\Pi_{n-1}^d$ with respect to the inner product $\la \cdot,\cdot \ra_{-s}$. Then $\CV_n^d(w_{-1})$ 
satisfies a decomposition 
$$
  \CV_n^d(w_{-1}) = (1-\|x\|^2) \CV_{n-2}^d(w_1)  \oplus \CH_n^d,
$$   
where $\CH_n^d$ is the space of spherical harmonics of degree $n$, and $\CV_n^d(w_{-2})$ satisfies a 
decomposition 
$$
  \CV_n^d(w_{-2})  = (1-\|x\|^2)^2 \CV_{n-4}^d(w_2)  \oplus (1-\|x\|^2) \CH_{n-2}^d \oplus \CH_n^d.
$$
For each of these two cases, an orthonormal basis can be given in terms of the Jacobi polynomials and 
spherical harmonics, and the basis resembles the basis of $\CV_n^d(w_\mu)$ for $\mu = -1$ and $\mu = -2$,
which is why we adopt the notation $\CV_n^d(w_{- s})$. The pattern of orthogonal decomposition, however, 
breaks down for $r >2$. Nevertheless, an orthonormal basis can still be defined for $\CV_n^d(w_{-s})$, which 
allows us to define an analogue of the near--best polynomial $S_{n,\eta}^{-s} f$. The result for approximation 
in the Sobolev space is as follows. 
 
\begin{thm}\label{nthApprox}
Let $r, s =1,2,\ldots$ and $r \ge s$. If $f \in W_p^{r}(\ball)$ with $r\ge s$ and
$1<p <\infty$. Then, for $n \ge s$, 
\begin{align*} 
 \| f -  S^{-s}_{n,\eta} f \|_{W_p^k(\ball)} \le  c n^{-r+k}  \|f\|_{W_p^r(\ball)},
 \quad k = 0, 1,\ldots, s, 
\end{align*}
where $S_{n, \eta}^{-s}f$  can be replaced by $S_{n}^{-s}f$ if $p=2$.
\end{thm}

This theorem is established in \cite{LX}, which contains further refinements of such estimates in Sobolev
spaces. The proof of this theorem, however, requires substantial work and uses a duality argument that 
requires $1 < p < \infty$. 

The estimate in the theorem can be used to obtain an error estimate for the Galerkin spectral method, which
looks for approximate solutions of a partial differential equations that are polynomials written in terms of 
orthogonal polynomials on the ball and their coefficients are determined by the Galerkin method. 
We refer to \cite{LX} for applications on a Helmholtz equation of second order and a biharmonic equation 
of fourth order on the unit ball. The method can also be applied to Poisson equations consider in 
\cite{ACH1, ACH2, AH}. 

These results raise the question of characterizing the best approximation by polynomials in Sobolev
spaces, which is closely related to simultaneous approximation traditionally studied in approximation
theory. But there are also distinct differences as the above discussion shows. We end this paper by formulating 
this problem in a more precise form. 

Let $\Omega$ be a domain in $\RR^d$ and $w$ be a weight function on $\Omega$. For $s =1,2\ldots$, 
and $f \in W_p^s(w, \Omega)$. Define 
$$
  E_n (f)_{W_p^s(w, \Omega)}: = \inf_{p_n \in \Pi_n^d} \|f - p_n \|_{W_p^s(w, \Omega)}. 
$$

\eject 
\noindent{\bf Problem 4.3.} Establish direct and (weak) converse estimates of $E_n (f)_{W_p^s(w, \Omega)}$. 
\medskip

In the case of $\Omega = \ball$ and $w (x) =1$, Theorem \ref{nthApprox} gives a direct estimate of
$E_n (f)_{W_p^s(w, \Omega)}$ for $f \in W_p^r(w, \Omega)$ with $r \ge s$. However, the estimate
is weaker than what is needed. A direct estimate should imply that $ E_n (f)_{W_p^s(w, \Omega)}$
goes to zero as $n \to \infty$ whenever $f \in W_p^s(w, \Omega)$. What this calls for is an appropriate 
$K$-functional, or a modulus of smoothness, for $f \in W_p^s(w, \Omega)$ that characterize the best
approximation $E_n (f)_{W_p^s(w, \Omega)}$.

\end{document}